\documentclass{ifacconf}

\usepackage{graphicx}      
\usepackage{natbib}        
\usepackage{amsmath,amssymb,bm,siunitx}
\usepackage{color}
\usepackage{fancyhdr}

\newcommand{\cov}[2]{\mathrm{cov}\left( #1, #2\right)}
\newcommand{\covs}[1]{\mathrm{cov}\left( #1 \right)}
\newcommand{\E}[1]{\mathrm{E}\left\{#1\right\}}
\newcommand{\mat}[1]{\left[\begin{matrix} #1 \end{matrix}\right]}
\newcommand{\norm}[1]{\left\lVert#1\right\rVert}
\DeclareMathOperator{\vecv}{vec}
\DeclareMathOperator{\diag}{diag}
\DeclareMathOperator{\ISOC}{ISOC}

\newtheorem{lemma}{Lemma}
\newtheorem{corollary}{Corollary}
\newtheorem{problem}{Problem}
\newtheorem{hypothesis}{Hypothesis}

\usepackage{ifthen}
\newboolean{showcomments}
\setboolean{showcomments}{false} 
\ifthenelse{\boolean{showcomments}}{ 
	\newcommand\td[1]{\textcolor{red}{\textbf{TODO:} #1}} 
	\newcommand\tdd[1]{} 
	\newcommand\comment[1]{\textcolor{blue}{\textbf{Comment:} #1}} 
	\newcommand\commentd[1]{} 
	\newcommand\frage[1]{\textcolor{orange}{\textbf{Rückfrage:} #1}} 
	\newcommand\fraged[1]{}
}{ 
	\newcommand\td[1]{}
	\newcommand\tdd[1]{} 
	\newcommand\comment[1]{} 
	\newcommand\commentd[1]{} 
	\newcommand\frage[1]{}
	\newcommand\fraged[1]{}
}

\begin{document}
\begin{frontmatter}

\title{Validation of Stochastic Optimal Control Models for Goal-Directed Human Movements on the Example of \\Human Driving Behavior}


\author[First]{Philipp Karg} 
\author[First]{Simon Stoll} 
\author[First]{Simon Rothfuß}
\author[First]{Sören Hohmann}

\address[First]{Institute of Control Systems (IRS), Karlsruhe Institute of Technology (KIT), Karlsruhe, Germany (e-mail: philipp.karg@kit.edu)}

\begin{abstract}                
Stochastic Optimal Control models represent the state-of-the-art in modeling goal-directed human movements. The linear-quadratic sensorimotor (LQS) model based on signal-dependent noise processes in state and output equation is the current main representative. With our newly introduced Inverse Stochastic Optimal Control algorithm building upon two bi-level optimizations, we can identify its unknown model parameters, namely cost function matrices and scaling parameters of the noise processes, for the first time. In this paper, we use this algorithm to identify the parameters of a deterministic linear-quadratic, a linear-quadratic Gaussian and a LQS model from human measurement data to compare the models' capability in describing goal-directed human movements. Human steering behavior in a simplified driving task shown to posses similar features as point-ot-point human hand reaching movements serves as our example movement. The results show that the identified LQS model outperforms the others with statistical significance. Particularly, the average human steering behavior is modeled significantly better by the LQS model. This validates the positive impact of signal-dependent noise processes on modeling human average behavior\comment{formal erstmal nur in der spezifischen driving task, da die movements in dieser aber vgl.bare kinematic features zeigen wie movements in point-to-point reaching tasks, kann durchaus step towards allgemeinerer conclusion erfolgen}.
\end{abstract}


\end{frontmatter}

\thispagestyle{fancy}
\pagestyle{fancy}
\fancyhf{}
\fancyhead[CO,CE]{\copyright 2023 IFAC. This paper has been accepted for publication and presentation at the IFAC World Congress 2023.}


\section{Introduction} \label{sec:intro}
\td{shorten introduction, first two sections of materials and methods, ggf. discussion and conclusion + shrink figures to achieve 7 pages in total}
Stochastic Optimal Control (SOC) models are the state-of-the-art approach to describe human movements to a single goal (\cite{Gallivan.2018}). Since they model basic optimality principles underlying human movements, like minimum intervention (\cite{Todorov.2002}), they are able to model the average behavior and the variability patterns of human limb movements. SOC models, and closed-loop optimization approaches in general, challenge traditional open-loop optimization models consisting of a separation between movement planning and execution (see e.g. \cite{Uno.1989,Flash.1985}). Moreover, due to their white-box character SOC models show higher generalizability and transferability between tasks than black-box models (alias behavioral cloning, see e.g. \cite{Ross.2010}). The current main SOC model is the linear-quadratic (LQ) sensorimotor (LQS) model, not least due to its computational tractability. The LQS model builds upon the well-known LQ Gaussian (LQG) framework. To this model, a control-dependent noise process is added to the state equation describing that faster movements are performed more inaccurately. Furthermore, a state-dependent noise process is added to the output equation describing the higher inaccuracy of the perception of faster movements. In order to validate SOC models in specific movement tasks, an identification method is needed to determine the unknown model parameters, namely relative weights of cost function candidates and scaling parameters of the noise processes, from human measurement data. Closely associated with the validation, such an identification method enables the analysis of optimality principles underlying human movements in a specific task as well (\cite{Franklin.2011,Todorov.2004}), i.e. preferred cost function candidates, influence of different noise processes on average and variability patterns as well as comparison to adapted models with omitted or added model parts. However, until very recently (\cite{Karg.2022}), such an identification method that identifies all unknown parameters of the LQS model from human measurement data was non-existent. The identified models and especially the conclusions drawn from the analysis of underlying human optimality principles are the crucial starting point for the design of human-machine systems that closely interact with the human in a cooperative, intuitive and safe manner: prediction, imitation and classification of human behavior is enabled as well as a model-based automation design.
\comment{Modellvalidierung bzw. Nachweis der Modellhypothese immer task-/movement-specific ebenso wie die Analyse der Optimalitätsprinzipien ("Beweise" für ein Modell sind meist generalisierender zu verstehen, sprich spiegelt Modell universelle Prinzipien menschl. Bewegungen wider, wie minimum intervention principle; datenbasiert über I(S)OC-Verfahren kann Modellhypothese dann nur task-/movement-specific belegt/analysiert werden); Experimente/Bewegungen müssen daher möglichst zahlreich, variierend und repräsentativ sein, letzteres schließt bspw. ein, dass komplexere Bewegungen aus diesen "Basisbewegungen" zusammengesetzt sind oder dass Bewegungen vgl.bare kinematic features zeigen, womit Modell ebenfalls wieder mit geeignet identifizierten Parametern geeignet sein sollte}

\cite{Kolekar.2018} motivate that human steering behavior is comparable to point-to-point reaching, widely considered as movement task in the neuroscientific community when human hand movements are studied. Since, in addition, a close cooperation between human and automation in steering a vehicle shows great potential to increase safety, for example in takeover scenarios (see \cite{Karg.2021}), human steering behavior in a simplified driving task acts as our example movement to investigate SOC models for goal-directed human movements. Describing human driving behavior via an optimal control framework traces back to (\cite{MacAdam.1981}) but only in recent years these optimal control models were extended to include neuroscientific findings (\cite{Kolekar.2018,Nash.2015,Markkula.2014,Sentouh.2009}). \cite{Nash.2015} propose the first stochastic human driver model that is based on the LQG framework. However, their model lacks signal-dependent noise processes which are promising to be considered (\cite{Kolekar.2018,Markkula.2014}). \cite{Kolekar.2018} make a first step to identify and validate the LQS model for human steering behavior but do not provide a rigor and generally applicable identification method that determines all model parameters. Moreover, they focus on describing the variability patterns of the human driver and hence, a comparison between the LQS, LQG and a deterministic LQ model regarding human average steering behavior is missing.

In this work, we want to deepen this line of research. First, we clarify that our steering task is comparable to a point-to-point reaching task by showing that the kinematic features of the movements in both tasks match. This strengthens the motivation done by \cite{Kolekar.2018} that steering is similar to reaching. Then, we use our newly introduced Inverse Stochastic Optimal Control (ISOC) algorithm (\cite{Karg.2022}) to identify all parameters of the LQS model from human measurement data and compare its modeling performance to reduced versions of it, i.e. to a LQG and LQ model. Here, our comparison and analysis is particularly focused on human average behavior and the influence of signal-dependent noise processes on it. This influence has theoretically been proven in our previous work (\cite{Karg.2022}).
\comment{hier bilden human measurement data die ground truth data (Input für ISOC-Algorithmus und ground truth für Evaluation, zu unterscheiden von reference einer menschlichen Bewegung im Experiment, in dem Mensch eine bestimmte Aufgabe/Task bekommt (Stichwort: "movement in driving task"), die dann gerade durch die reference beschrieben wird; GT data kann i.Allg. aus Simulation oder aufgenommenen Messdaten kommen); da hier stets ground truth = human measurement gilt, hier nur human measurement data nutzen; in der human measurement data werden dann nicht alle Größen des Modells gemessen/measured (dies aber nicht als partially measurable o.Ä. beschreiben); partially/fully observable bezieht sich dann ganz klass. nur auf das LQS optimal control problem zur Beschreibung der human perception}

\tdd{
	\begin{itemize}
		\item human movement modeling (see MindMap, SOC and closed-loop = SoA \cite{Gallivan.2018,Todorov.2002}, more general than black-box approaches (alias behavioral cloning) \cite{Ross.2010} and separation between movement planning and execution (alias open-loop approaches) \cite{Flash.1985,Uno.1989} challenged in the last years)
		\item suitable identification approaches (for all model parameters) missing until recently \cite{Karg.2022}
		\item necessary to finally validate and analyze model hypothesis (done task-/movement-specific, see \cite{Franklin.2011,Todorov.2004}): analyze best relative weights of cost function candidates (optimality principles underlying human movements), analyze influence of the different noise processes and their parameters, comparison to adapted models (omission or addition of model parts)
		\item identified models and drawn conclusions beneficial for the design of human-machine systems that closely interact with the human in a cooperative, intuitive and safe manner: prediction, imitation and classification of human behavior + model-based automation design
		\item validation and analysis of the LQS model hypothesis on the example of human driving behavior (steering task), comparison to reduced LQG and deterministic LQ model
		\item human driver models based on OC framework common \cite{MacAdam.1981}; in recent years, they were extended to include neuroscientific findings \cite{Sentouh.2009,Markkula.2014,Nash.2015,Kolekar.2018}; \cite{Nash.2015} LQG framework which is novel compared to traditional uncertainty approaches in modeling driving behavior \cite{Gray.2013} but lack signal-dependent noise processes which are very promising to consider \cite{Markkula.2014,Kolekar.2018}; \cite{Kolekar.2018} first to identify/validate LQS model in this case but no rigor/general identification approach that identifies all model parameters, focus on variability, comparison between LQ, LQG and LQS missing; furthermore, they provide a first motivation for steering as reaching; in this paper, we want to deepen this line of research and provide a validation of SOC models for goal-directed human movements on the example of a simplified steering task (with our newly developed ISOC algorithm)
		\item contributions: full SOC models are identified/trained on human movement data and compared between each other and to deterministic model, focus on comparison/analysis regarding the mean (and the influence of signal-dependent noise parameters on it, influence proved in \cite{Karg.2022}), steering similar to reaching (shown by analysis of kinematic features, more profound than motivation in \cite{Kolekar.2018})
	\end{itemize}
}


\section{Materials and Methods} \label{sec:materials_methods}
In the following section the LQS model and our ISOC algorithm (\cite{Karg.2022}) are introduced. Afterwards, the experimental setup is explained. 
\tdd{up to and including page 3}

\subsection{Stochastic Optimal Control Models for Goal-Directed Human Movements} \label{subsec:soc}
The human biomechanics and the system the human interacts with, e.g. a vehicle with steering wheel, is modeled via the linear time-discrete system dynamics
\begin{align}
	\bm{x}_{t+1} &= \bm{A} \bm{x}_{t} + \bm{B} \bm{u}_{t} + \bm{\Sigma}^{\bm{\xi}}\bm{\alpha}_{t} + \sum_{i=1}^{c} \varepsilon_t^{(i)} \bm{C}_i \bm{u}_t \label{eq:stateEquationLQS} \\
	\bm{y}_{t} &= \bm{H} \bm{x}_{t} + \bm{\Sigma}^{\bm{\omega}}\bm{\beta}_{t}  + \sum_{i=1}^{d} \epsilon_t^{(i)} \bm{D}_i \bm{x}_t \label{eq:outputEquationLQS},
\end{align}
where $\bm{x}\in \mathbb{R}^{n}$ denotes the system state, $\bm{u}\in \mathbb{R}^m$ the control variable (typically neural activation of muscles), $\bm{y}\in \mathbb{R}^r$ the human observation and $\bm{A}, \bm{B}, \bm{H}$ matrices of appropriate dimension that may depend on time. When the additional time index $t$ is used, e.g. $\bm{x}_t$, the stochastic process of the corresponding variable is denoted. Furthermore, $\bm{\alpha}_t \in \mathbb{R}^p$ and $\bm{\varepsilon}_t = \mat{\varepsilon_t^{(1)} & \dots & \varepsilon_t^{(c)}}^{\top}$ are standard white Gaussian noise processes that constitute the additive and control-dependent noise process of state equation~\eqref{eq:stateEquationLQS} with their scaling matrices $\bm{\Sigma}^{\bm{\xi}}$ and $\bm{C}_i = \sigma^{\bm{u}}_i \bm{B}\bm{F}_i$ of appropriate dimension. In \eqref{eq:outputEquationLQS}, $\bm{\beta}_t \in \mathbb{R}^q$ and $\bm{\epsilon}_t = \mat{\epsilon_t^{(1)} & \dots & \epsilon_t^{(d)}}^{\top}$ are standard white Gaussian noise processes forming the additive and state-dependent noise process of the output equation with their corresponding scaling matrices $\bm{\Sigma}^{\bm{\omega}}$ and $\bm{D}_i = \sigma^{\bm{x}}_i \bm{H}\bm{G}_i$ of appropriate dimension. The stochastic processes $\bm{\alpha}_t$, $\bm{\varepsilon}_t$, $\bm{\beta}_t$ and $\bm{\epsilon}_t$ are independent to each other and to $\bm{x}_t$. The goal of the human, e.g. achievement of a desired steering angle under minimal effort, is modeled by the performance criterion 
\begin{align}
	J = \E{\bm{x}^{\top}_N \bm{Q}_N \bm{x}_N + \sum_{t=0}^{N-1} \bm{x}^{\top}_t \bm{Q} \bm{x}_t + \bm{u}^{\top}_t \bm{R} \bm{u}_t}, \label{eq:costLQS}
\end{align} 
where $\bm{Q}_N$, $\bm{Q}$ and $\bm{R}$ are symmetric and of appropriate dimension. $\bm{Q}_N$, $\bm{Q}$ are positive semi-definite and $\bm{R}$ is positive definite. It is assumed that the human solves the LQS control problem consisting of the minimization of \eqref{eq:costLQS} with respect to \eqref{eq:stateEquationLQS} and \eqref{eq:outputEquationLQS} by an admissible control strategy. An approximate (suboptimal) solution to the LQS control problem is given by \cite{Todorov.2005} with $\bm{u}_t = -\bm{L}_t \hat{\bm{x}}_t$ and $\hat{\bm{x}}_{t+1} = \bm{A}\hat{\bm{x}}_t + \bm{B}\bm{u}_t + \bm{K}_t (\bm{y}_t - \bm{H}\hat{\bm{x}}_t) + \bm{\eta}_t$, where $\bm{L}_t$ and $\bm{K}_t$ are computed by iterating between 
\begin{align}
	\bm{L}_{t} &= \Big(\bm{R} + \bm{B}^{\top} \bm{Z}_{t+1}^{\bm{x}} \bm{B} \nonumber\\
	&\hphantom{=} + \sum_{i} \bm{C}_i^{\top} \left( \bm{Z}_{t+1}^{\bm{x}} + \bm{Z}_{t+1}^{\bm{e}} \right) \bm{C}_i \Big)^{-1} \bm{B}^{\top} \bm{Z}_{t+1}^{\bm{x}} \bm{A} \label{eq:controlLawLQS} 
\end{align}
and
\begin{align}
	\bm{K}_{t} &= \bm{A} \bm{P}_{t}^{\bm{e}} \bm{H}^{\top} \Big( \bm{H} \bm{P}_{t}^{\bm{e}} \bm{H}^{\top} + \bm{\Omega}^{\bm{\omega}} \nonumber \\
	&\hphantom{=} + \sum_{i} \bm{D}_i \left( \bm{P}_t^{\bm{e}} + \bm{P}_t^{\hat{\bm{x}}} + \bm{P}_t^{\hat{\bm{x}}\bm{e}} + \bm{P}_t^{\bm{e}\hat{\bm{x}}} \right) \bm{D}_i^{\top} \Big)^{-1} \label{eq:estimatorLQS}.
\end{align}
In \eqref{eq:controlLawLQS} and \eqref{eq:estimatorLQS}, $\bm{Z}_{t}^{\bm{e}}$, $\bm{Z}_{t}^{\bm{x}}$, $\bm{P}_t^{\bm{e}}$, $\bm{P}_t^{\hat{\bm{x}}}$, $\bm{P}_t^{\bm{e}\hat{\bm{x}}}$ and $\bm{P}_t^{\hat{\bm{x}}\bm{e}}$ follow from recursive equations (see \cite{Todorov.2005}).

Suppose the cost function matrices $\bm{Q}_N$, $\bm{Q}$, $\bm{R}$ and the noise scaling parameters $\bm{\Sigma}^{\bm{\xi}}$, $\sigma_i^{\bm{u}}$, $\bm{\Sigma}^{\bm{\omega}}$, $\sigma_i^{\bm{x}}$ are known, the average human behavior (mean $\E{\bm{x}_t}$) and the human variability pattern (covariance $\bm{\Omega}_t^{\bm{x}} = \cov{\bm{x}_t}{\bm{x}_t} = \covs{\bm{x}_t}$) predicted by the LQS model are given by Lemma~\ref{lemma:solutionLQS}.
\begin{lemma} \label{lemma:solutionLQS} 
	Let the LQS control problem be defined by the minimization of \eqref{eq:costLQS} with respect to \eqref{eq:stateEquationLQS} and \eqref{eq:outputEquationLQS}. Let the (approximate) solution be given by $\bm{L}_t$~\eqref{eq:controlLawLQS} and $\bm{K}_t$~\eqref{eq:estimatorLQS}. Then, the mean $\E{\bm{x}_t}$ and covariance $\bm{\Omega}_t^{\bm{x}}$ of $\bm{x}_t$ are 
	\begin{align}
		\mat{\E{\bm{x}_{t+1}} \\ \E{\hat{\bm{x}}_{t+1}}} &= \bm{\mathcal{A}}_t \mat{\E{\bm{x}_{t}} \\ \E{\hat{\bm{x}}_{t}}} \label{eq:meanLQS}, \\
		\mat{\bm{\Omega}_{t+1}^{\bm{x}} & \bm{\Omega}_{t+1}^{\bm{x}\hat{\bm{x}}} \\ \bm{\Omega}_{t+1}^{\hat{\bm{x}}\bm{x}} & \bm{\Omega}_{t+1}^{\hat{\bm{x}}}} &= \bm{\mathcal{A}}_t \mat{\bm{\Omega}_{t}^{\bm{x}} & \bm{\Omega}_{t}^{\bm{x}\hat{\bm{x}}} \\ \bm{\Omega}_{t}^{\hat{\bm{x}}\bm{x}} & \bm{\Omega}_{t}^{\hat{\bm{x}}}} \bm{\mathcal{A}}^\top_{t} \nonumber \\		
		&\hphantom{=} + \mat{\bm{\Omega}^{\bm{\xi}} & \bm{0} \\ \bm{0} & \bm{\Omega}^{\bm{\eta}} \! + \! \bm{K}_{t} \bm{\Omega}^{\bm{\omega}} \bm{K}^{\top}_{t}} + \mat{\bar{\bm{\Omega}}_t^{\hat{\bm{x}}} & \bm{0} \\ \bm{0} & \bar{\bm{\Omega}}_t^{\bm{x}}} \label{eq:covarianceLQS}
	\end{align}
	with $\bar{\bm{\Omega}}_t^{\hat{\bm{x}}} = \sum_{i} \bm{C}_i \bm{L}_t \left( \bm{\Omega}_{t}^{\hat{\bm{x}}} + \E{\hat{\bm{x}}_t}\E{\hat{\bm{x}}_t}^{\top} \right) \bm{L}_t^{\top} \bm{C}_i^{\top}$, $\bar{\bm{\Omega}}_t^{\bm{x}} = \sum_{i} \bm{K}_t \bm{D}_i \left( \bm{\Omega}_{t}^{\bm{x}} + \E{\bm{x}_t}\E{\bm{x}_t}^{\top} \right) \bm{D}_i^{\top} \bm{K}_t^{\top}$, 
	\begin{align}
		\bm{\mathcal{A}}_{t} = \mat{\bm{A} & -\bm{B} \bm{L}_{t} \\ \bm{K}_{t} \bm{H} & \bm{A} - \bm{K}_{t} \bm{H} - \bm{B} \bm{L}_{t}} \label{eq:mathcalA}
	\end{align}
	and initial values $\E{\hat{\bm{x}}_0}=\hat{\bm{x}}_0=\E{\bm{x}_0}$ and
	\begin{align}
		\mat{\bm{\Omega}_{0}^{\bm{x}} & \bm{\Omega}_{0}^{\bm{x}\hat{\bm{x}}} \\ \bm{\Omega}_{0}^{\hat{\bm{x}}\bm{x}} & \bm{\Omega}_{0}^{\hat{\bm{x}}}} = \mat{\bm{\Omega}_{0}^{\bm{x}} & \bm{0} \\ \bm{0} & \bm{0}}.
	\end{align}
\end{lemma}
\begin{pf}
	See \cite[Lemma~2]{Karg.2022}.
\end{pf}
\begin{corollary} \label{corollary:meanLQS}
	The mean $\E{\bm{x}_t}$ of the closed-loop LQS system depends on the noise scaling parameters $\sigma_i^{\bm{u}}$ in the fully observable case ($\hat{\bm{x}}_t = \bm{x}_t$) and on $\sigma_i^{\bm{u}}$, $\sigma_i^{\bm{x}}$, $\bm{\Sigma}^{\bm{\xi}}$ and $\bm{\Sigma}^{\bm{\omega}}$ in the partially observable case.
\end{corollary}
\begin{pf}
	From \eqref{eq:meanLQS} and \eqref{eq:mathcalA} the dependency of $\E{\bm{x}_t}$ on $\bm{L}_t$~\eqref{eq:controlLawLQS} follows. In the fully observable case, $\bm{Z}_t^{\bm{e}}$ drops out in \eqref{eq:controlLawLQS} (see \cite{Todorov.2005}) but $\bm{L}_t$ still depends on $\bm{C}_i$ and thus $\sigma_i^{\bm{u}}$. In the partially observable case, the additional dependencies follow from the dependency of $\bm{L}_t$ on $\bm{Z}_t^{\bm{e}}$ which in turn depends on $\bm{K}_t$ (see \cite{Todorov.2005}).
\end{pf}

Corollary~\ref{corollary:meanLQS} shows that the human average behavior $\E{\bm{x}_t}$ predicted by the LQS model depends not only on the cost function matrices but also on the scaling parameters of the noise processes. Hence, they could act as a crucial part of the LQS model improving the performance in describing human average behavior.

\tdd{
\begin{itemize}
	\item LQS model
	\item repeat lemma (without proof), shows influence of signal-dependent noise processes on mean (average behavior)
	\item LQG and LQ follow directly by simplifications (alias reduced versions of the LQS model)
\end{itemize}
}

\subsection{Inverse Stochastic Optimal Control} \label{subsec:isoc}
Regarding the LQS model introduced in the previous section, the cost function matrices $\bm{Q}_N$, $\bm{Q}$, $\bm{R}$ and the noise (scaling) parameters $\bm{\Sigma}^{\bm{\xi}}$, $\sigma_i^{\bm{u}}$, $\bm{\Sigma}^{\bm{\omega}}$, $\sigma_i^{\bm{x}}$ are unknown and need to be identified when looking at human movements in practice. To simplify notation, we combine the cost function parameters in a vector $\bm{s} \in \mathbb{R}^{S}$ and the noise parameters in a vector $\bm{\sigma} \in \mathbb{R}^{\Sigma}$ (see \cite{Karg.2022}). Now, we assume that our measurement data consists of a set $\{\bm{M}\bm{x}^{\ast,(k)}_{t}\}$ ($k \in \{1,\dots,K\}$) of time-discrete trajectories $\bm{M}\bm{x}^{\ast,(k)}_{t}$ where $\bm{x}^{\ast,(k)}_{t}$ is a realization of the stochastic process $\bm{x}_t^{\ast}$. It results in the closed-loop LQS system with the unknown cost function $\bm{s}^{\ast}$ and noise parameters $\bm{\sigma}^{\ast}$ of the human. With $\bm{M}\in \mathbb{R}^{\bar{n} \times n}$, we consider that not all system states are measured, i.e. $\bm{M}$ follows from the identity matrix by deleting rows corresponding to not measured states. The identification problem is finally formalized in Problem~\ref{problem:ISOC}\footnote{The problem extends typical Inverse Optimal Control or Inverse Reinforcement Learning problems by additionally searching for scaling parameters of noise processes.}.
\begin{problem} \label{problem:ISOC}
	Let estimates $\hat{\bm{m}}_t \approx \E{\bm{M}\bm{x}_t^{\ast}}$ and $\hat{\bm{\Omega}}^{\bm{x}^{\ast}}_{t} \approx \bm{M}\bm{\Omega}^{\bm{x}^{\ast}}_{t}\bm{M}^{\top}$ computed from the measurement data be given. Find parameters $\tilde{\bm{s}}$ and $\tilde{\bm{\sigma}}$ such that they lead to $\tilde{\bm{x}}_t$ in the closed-loop LQS system with $\E{\bm{M}\tilde{\bm{x}}_t}=\E{\bm{M}\bm{x}^{\ast}_t}$ and $\bm{M}\bm{\Omega}^{\tilde{\bm{x}}}_{t}\bm{M}^{\top}=\bm{M}\bm{\Omega}^{\bm{x}^{\ast}}_{t}\bm{M}^{\top}$.
\end{problem}

\begin{figure}[t]
	\centering
	\includegraphics[width=3.5in]{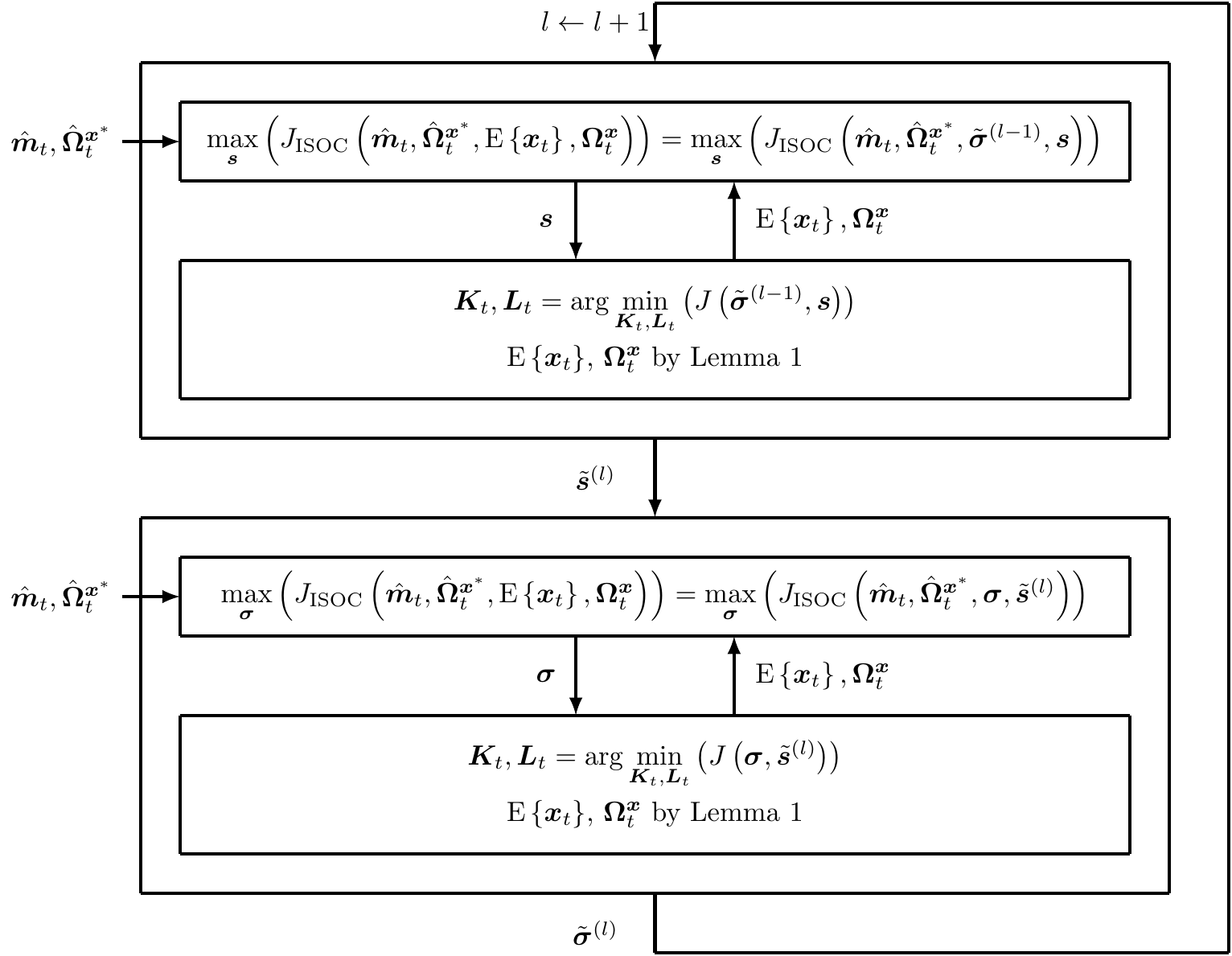}
	\caption{Bi-level-based ISOC algorithm.}
	\label{fig:isoc}
\end{figure}

In order to solve Problem~\ref{problem:ISOC}, we have introduced an ISOC algorithm based on two bi-level optimizations in our previous work (\cite{Karg.2022}). We formulate a direct optimization problem by defining a performance criterion $J_{\ISOC}$ describing how well the mean $\E{\bm{x}_t}$ and covariance $\bm{\Omega}_t^{\bm{x}}$ predicted by the LQS model with parameters $\bm{s}$ and $\bm{\sigma}$ match the measurement data $\hat{\bm{m}}_t$ and $\hat{\bm{\Omega}}_t^{\bm{x}^{\ast}}$. The objective function $J_{\ISOC}$ is chosen based on the Variance Accounted For (VAF) metric:
\begin{align}
	J_{\ISOC} = \frac{\bm{w}_{\text{m}}^{\top} \bm{m}^{\text{VAF}} + \bm{w}_{\text{v}}^{\top} \vecv\left(\bm{\Omega}^{\text{VAF}} \right)}{\norm{\bm{w}_{\text{m}}}_1 + \norm{\bm{w}_{\text{v}}}_1} \label{eq:J_ISOC},
\end{align}
where for the elements $m^{\text{VAF}}_i$ and $\Omega^{\text{VAF}}_{ij}$ 
\begin{align}
	m^{\text{VAF}}_i &= \left(1 - \frac{\sum_{t=0}^{N} \left(\left(\E{\bm{M}\bm{x}_t}\right)_i - \hat{m}_{i,t}\right)^2}{\sum_{t=0}^{N} \left(\hat{m}_{i,t} - \frac{1}{N+1}\sum_{t}\hat{m}_{i,t}\right)^2}\right) \label{eq:VAFEW}, \\
	\Omega^{\text{VAF}}_{ij} &= \left(1 - \frac{\sum_{t=0}^{N} \left(\left(\bm{M}\bm{\Omega}^{\bm{x}}_{t}\bm{M}^{\mathrm{T}}\right)_{ij} - \hat{\Omega}_{ij,t}^{\bm{x}^{\ast}}\right)^2}{\sum_{t=0}^{N} \left(\hat{\Omega}_{ij,t}^{\bm{x}^{\ast}} - \frac{1}{N+1}\sum_{t}\hat{\Omega}_{ij,t}^{\bm{x}^{\ast}}\right)^2}\right) \label{eq:VAFCOV}
\end{align}
holds. In \eqref{eq:J_ISOC}, $\bm{w}_{\text{m}} \in \mathbb{R}^{\bar{n}}$ and $\bm{w}_{\text{v}} \in \mathbb{R}^{\bar{n}\bar{n}}$ denote weighting vectors. Due to the VAF metric $J_{\ISOC} \in (-\infty,1]$ holds with $1$ corresponding to a perfect fit between predicted and measured data. Therefore, the direct ISOC optimization problem consists of maximizing $J_{\ISOC}$: $\max \limits_{\bm{s},\bm{\sigma}} \left( J_{\ISOC}\left( \hat{\bm{m}}_t, \hat{\bm{\Omega}}^{\bm{x}^{\ast}}_t, \E{\bm{x}_t}, \bm{\Omega}^{\bm{x}}_t \right) \right)$. The evaluation of $J_{\ISOC}$ with specific values for $\bm{s}$ and $\bm{\sigma}$ depends on calculating the model predictions $\E{\bm{x}_t}$ and $\bm{\Omega}_t^{\bm{x}}$ which in turn depends on solving the forward optimal control problem, i.e. computing $\bm{K}_t$ and $\bm{L}_t$. Hence, maximizing $J_{\ISOC}$ leads to a bi-level optimization problem where the upper level optimization consists of maximizing $J_{\ISOC}$ and the lower level of calculating  $\bm{K}_t$, $\bm{L}_t$, $\E{\bm{x}_t}$ and $\bm{\Omega}_t^{\bm{x}}$. In order to keep this optimization computationally tractable, we apply an alternating optimization approach by computing the best estimate $\tilde{\bm{s}}^{(l)}$ for a given estimate $\tilde{\bm{\sigma}}^{(l-1)}$ in the first step of one iteration $l$ and in the second vice versa. This leads to the iteration between two bi-level optimizations shown in Fig.~\ref{fig:isoc}. Furthermore, in the lower levels of both bi-level optimizations our newly introduced Lemma~\ref{lemma:solutionLQS} is utilized to simplify the computation of $\E{\bm{x}_t}$ and $\bm{\Omega}_t^{\bm{x}}$ by applying recursive calculations only. Finally, for solving the upper level optimizations the parameters $\bm{s}$ and $\bm{\sigma}$ are divided into parameter sets where each parameter is supposed to be in at least one set and where the parameters in one set are assumed to have a mutually correlated influence on $J_{\ISOC}$. Then, grid searches with fixed grid size are performed on these parameter sets iteratively until no further improvement with the current grid size is achieved. Then, the grid size is shrinked. Further details regarding the ISOC algorithm can be found in (\cite{Karg.2022}).

\tdd{
	\begin{itemize}
		\item parameters that need to be identified, measured ground truth data
		\item introduce direct optimization objective and VAF metric
		\item algorithm from \cite{Karg.2022} (in improved version)
		\item figure, but slightly more abstract
	\end{itemize}
}

\subsection{Human Driving Experiment} \label{subsec:driving_experiment}
\begin{figure}[t]
	\centering
	\includegraphics[width=2.5in]{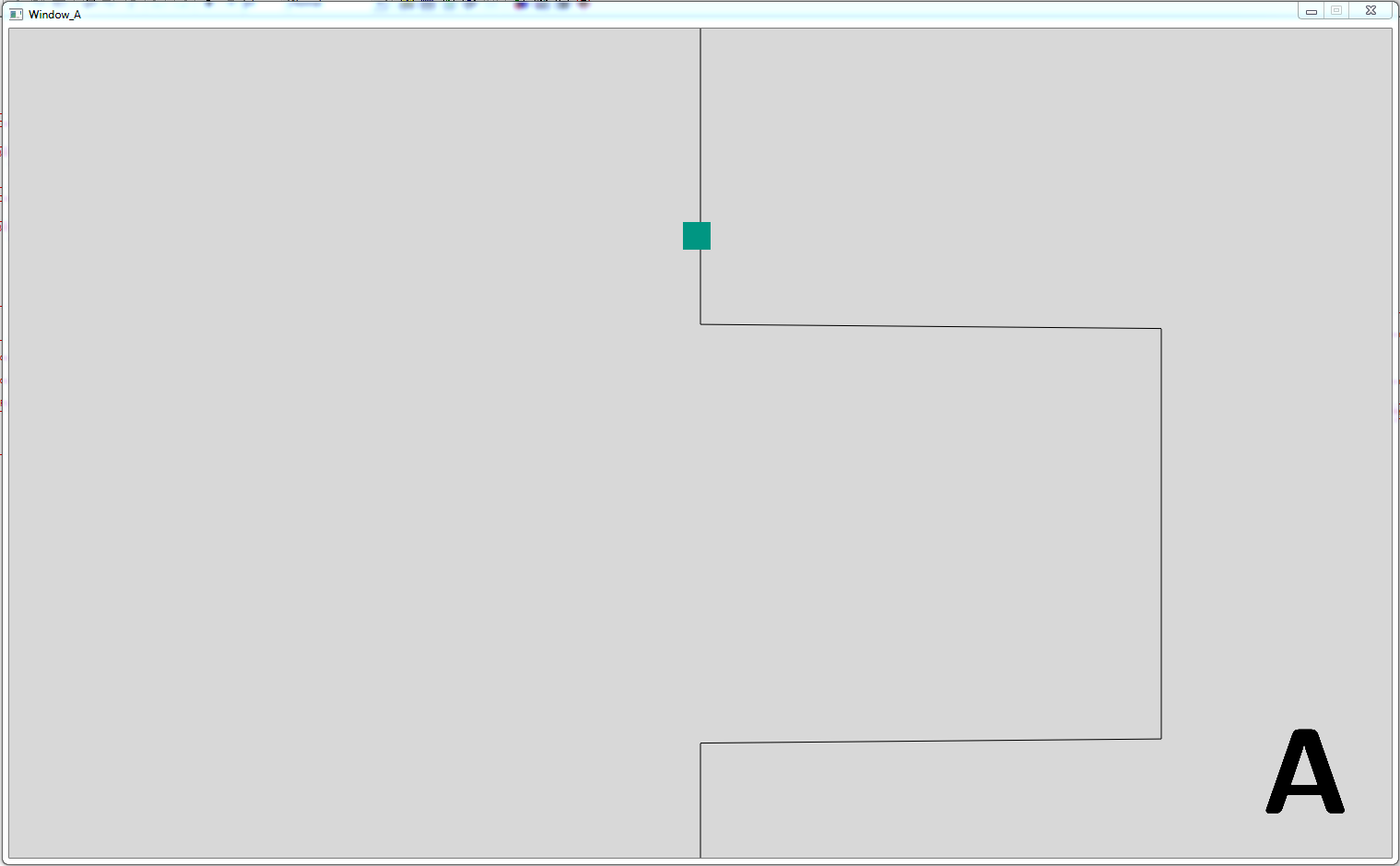}
	\caption{Screen visualization of the simplified steering task. The study participants interacted with an active steering wheel to change the horizontal position of the squared marker and were asked to follow the piecewise constant reference trajectory as best as possible.}
	\label{fig:screenVisualization}
\end{figure}

\tdd{clarify introduction to drivers and familiarization phase, add figure of representations on screen and connect it to the explanations}
The human measurement data was gathered by a study where 14 participants (12 male and 2 female subjects) aged between 19 and 29 (24.9 in average) performed a simplified steering task. The participants interacted with an active steering wheel which was equipped with an incremental encoder (40000 steps per full rotation) to measure the steering angle $\varphi$ with a sampling frequency of $\SI{100}{Hz}$. Based on the steering angle measurements torque was applied by the steering wheel such that it mimicked the dynamics of a spring-damper system. With the steering wheel the participants were able to move a marker (square) on a screen on a horizontal line with fixed vertical position (cf. Figure~\ref{fig:screenVisualization}). The horizontal position of the marker was a linear mapping of the steering angle.
On the display, a reference trajectory for the marker (cf. Figure~\ref{fig:screenVisualization}) was shown and the participants were asked to follow it as best as possible. The reference trajectory moved downwards with constant velocity. The reference trajectory shown for the marker is directly connected to a reference trajectory $\varphi_{\text{r},t}$ for the steering angle. The trial of one participant started with a familiarization phase of $\SI{210}{s}$ consisting of different segments with a steering angle reference with constant curvature where the curvature varied between those segments. Afterwards, the reference consisted of 14 repetitions of the same piecewise constant pattern: $\varphi_{\text{r}} = 0$, $\varphi_{\text{r}} = +\frac{2}{3}\pi$, $\varphi_{\text{r}} = 0$ and $\varphi_{\text{r}} = -\frac{2}{3}\pi$, each for a duration of $100$ time steps. 

Regarding the data preparation, we first applied a cubic spline smoothing on the measured steering angle to reduce measurement noise. Afterwards, we calculate the steering angle velocity $\dot{\varphi}_t$ via numerical differentiation. Since the measured data consists of 14 repetitions of the $+\frac{2}{3}\pi$-step as well as 14 repetitions of the $-\frac{2}{3}\pi$-step, we calculate estimates of human average behavior ($\hat{\bm{m}}_t$) and human variability pattern ($\hat{\bm{\Omega}}_t^{\bm{x}^\ast}$) for each of these movements from the corresponding 14 repetitions. Hereto, we determine the starting time for each single movement repetition as the last point in time where $\varphi \approx 0$ and $\dot{\varphi} < \SI{0.1}{s^{-1}}$ holds. In some segments with $\varphi_{\text{r}}=0$ the steering angle velocity is always above this threshold. Then, we increase the threshold in $\SI{0.01}{s^{-1}}$-steps until a maximum value of $\SI{0.4}{s^{-1}}$. If even this maximally relaxed condition is not fulfilled, we remove this single movement repetition from the further analysis since the human clearly starts not from an equilibrium condition in this case. Finally, the remaining repetitions are used to calculate the mean ($\hat{\bm{m}}_t$) and covariance estimate ($\hat{\bm{\Omega}}_t^{\bm{x}^\ast}$) of $\varphi_t$ and $\dot{\varphi}_t$. The movement duration is defined as the averaged movement duration of those repetitions which fulfill the condition for a movement ending ($\varphi \approx \pm \frac{2}{3}\pi$ and $\dot{\varphi} < \SI{0.1}{s^{-1}}$). This data preparation is done for every subject.

\begin{figure}[t]
	\centering
	\includegraphics[width=2.5in]{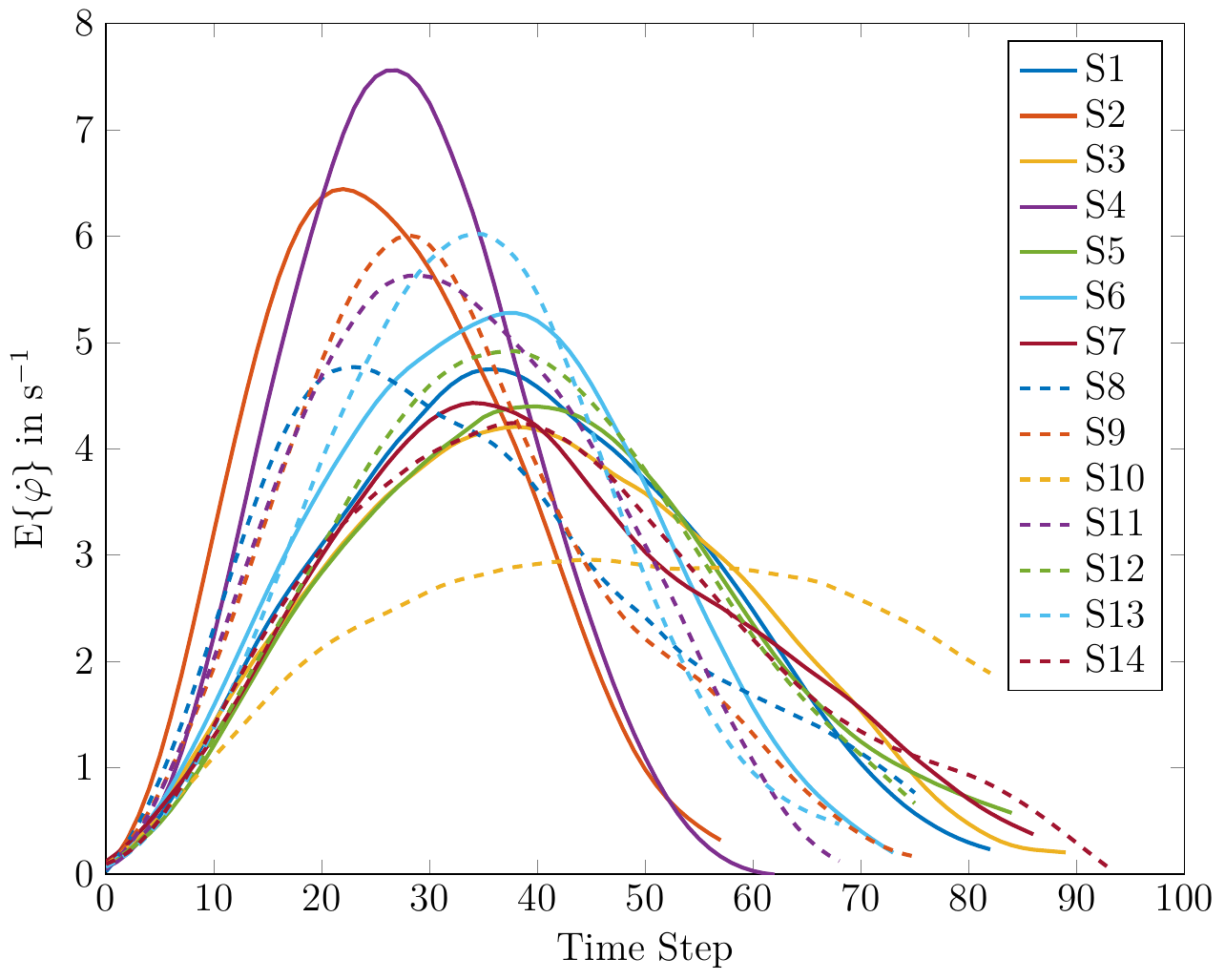}
	\caption{Steering angle velocity profiles of all 14 subjects (S1-S14) for the $+\frac{2}{3}\pi$-step.}
	\label{fig:velocityProfiles}
\end{figure}

Fig.~\ref{fig:velocityProfiles} shows the steering angle velocity profiles for the $+\frac{2}{3}\pi$-step achieved by the averaging procedure described before. Remarkably, the profiles possess very similar kinematic features as point-to-point reaching movements: they are all single-peaked and predominantly bell-shaped and nearly symmetric. The majority of the profiles is furthermore slightly left-skewed which indicates that the subjects tend to perform fast movements (\cite{Engelbrecht.2001}). Moreover, the maximum to mean ratios of $\E{\dot{\varphi}}$ range from $1.78$ to $2.26$ (when S10 is not considered due to its outlier behavior). All these features were also observed for point-to-point human hand movements (\cite{Engelbrecht.2001}). This conclusion shows on one side the validity of our measured data and on the other that we can make a first step towards a data-driven validation of SOC models for general goal-directed human movements by applying our ISOC algorithm. 

Before results of this application are presented in the next section, the system equations \eqref{eq:stateEquationLQS} and \eqref{eq:outputEquationLQS} need to be set up for the driving experiment. The active steering wheel is modeled by the spring-damper dynamics
\begin{equation}
	\Theta\ddot{\varphi} = -c\varphi -d\dot{\varphi} + M_{\text{h}} \label{eq:springDamper},
\end{equation}
where $M_{\text{h}}$ is the steering torque applied by the human. It results from a second-order linear filter with time constants $\tau_1$ and $\tau_2$ describing the low-pass characteristics of the muscle dynamics (\cite{Winter.1990}):
\begin{equation}
	\tau_1 \dot{g} + g = u, \quad \tau_2 \dot{M_\text{h}} + M_\text{h} = g \label{eq:muscleFilter},
\end{equation}
where $u$ is the neural activation and $g$ the muscle excitation. We define the system state as $\bm{x} = \mat{\varphi & \dot{\varphi} & M_{\text{h}} & g & \varphi_{\text{r}}}^{\top}$. The augmentation of it with the target steering angle $\varphi_{\text{r}}$ (constant dynamics) is a mathematical reformulation to consider target steering angles $\varphi_{\text{r}} \neq 0$ but still constitute the cost function in the form of \eqref{eq:costLQS}:
\begin{equation}
	J = \E{ s_1 (x_{N,1} - x_{N,5})^2 + s_2 x_{N,2}^2 + s_3 x_{N,3}^2 + \sum_{t=0}^{N-1} s_4 u_t^2} \label{eq:drivingExpJ},
\end{equation}
from which the definition of the cost function vector $\bm{s}$ follows. Discretizing \eqref{eq:springDamper} and \eqref{eq:muscleFilter} and choosing the numerical values as $\Theta=\SI{0.056}{kg m^2}$, $c=\SI{1.146}{Nm}$, $d=\SI{0.859}{Nms}$ and $\tau_1=\tau_2=\SI{0.04}{s}$ yields $\bm{A}$ and $\bm{B}$. Furthermore, $\bm{H} = \mat{\bm{I}_{3\times3} & \bm{0}_{3\times2}}$ is assumed for the human perception since $\varphi_{\text{r}}$ only results from the above mentioned mathematical reformulation and $g$ cannot be clearly associated with the human sensory apparatus. Due to lack of further knowledge about the additive noise processes and to keep the number of noise parameters tractable, independent stochastic processes $\alpha_{i,t}$ and $\beta_{i,t}$ are modeled for every state and output: $\bm{\Sigma}^{\bm{\xi}} = \diag(\mat{\sigma_1 & \sigma_2 & \sigma_3 & \sigma_4 & 0})$ and $\bm{\Sigma}^{\bm{\omega}} = \diag(\mat{\sigma_5 & \sigma_6 & \sigma_7})$. Finally, the signal-dependent noise processes are modeled by $\bm{C} = \sigma_8\bm{B}$ (single-input system), $\bm{D}_1 = \sigma_9 \bm{H} \diag(\mat{1 & 0 & 0 & 0 & 0})$, $\bm{D}_2 = \sigma_{10} \bm{H} \diag(\mat{0 & 1 & 0 & 0 & 0})$ and $\bm{D}_3 = \sigma_{11} \bm{H} \diag(\mat{0 & 0 & 1 & 0 & 0})$ which completes the definition of the noise parameter vector $\bm{\sigma}$. Again, due to lack of further knowledge about the state-dependent noise processes, the $\bm{G}_i$ are chosen such that every perceived state is influenced by an independent noise process scaled by the corresponding state. In case of the LQG and LQ model, the signal-dependent and all noise parameters drop out, respectively. 

For applying our ISOC algorithm, we first divide the measurement data into training and validation data. The parameters $\bm{s}$ and $\bm{\sigma}$ of the LQS, LQG and LQ model are identified for each subject individually with $\hat{\bm{m}}_t$ and $\hat{\bm{\Omega}}_t^{\bm{x}^{\ast}}$ calculated from the repetitions of the $+\frac{2}{3}\pi$-step. Here, we set $\bm{M}=\mat{\bm{I}_{2\times2} & \bm{0}_{2\times3}}$ since $\varphi$ and $\dot{\varphi}$ are measured. The repetitions of the $-\frac{2}{3}\pi$-step serve as validation data, i.e. the corresponding values for $\hat{\bm{m}}_t$ and $\hat{\bm{\Omega}}_t^{\bm{x}^{\ast}}$ are the ground truth data for the LQS, LQG and LQ model predictions computed with $\tilde{\bm{s}}$ and $\tilde{\bm{\sigma}}$. The weighting vectors for the cost function bi-level optimization are chosen as $\bm{w}_{\text{m}} = \mat{0.9 & 0.9}^{\top}$ and $\bm{w}_{\text{v}} = \mat{0.1 & 0}^{\top}$ in case of the LQS and LQG model and as $\bm{w}_{\text{m}} = \mat{1 & 1}^{\top}$ and $\bm{w}_{\text{v}} = \mat{0 & 0}^{\top}$ for the LQ model. For the noise parameter bi-level optimization in case of the LQS and LQG model, $\bm{w}_{\text{m}} = \mat{0.1 & 0.1}^{\top}$ and $\bm{w}_{\text{v}} = \mat{0.9 & 0}^{\top}$ holds. 

\tdd{
	\begin{itemize}
		\item experimental setup, explain training and validation data: 14 participants (track type $\beta$, steering system B), 14 trials with one $+120^{\circ}$ and one $-120^{\circ}$ step in target steering angle per participant; no additional trials from movement to zero angle equilibrium due to its equilibrium character; 14 $+120^{\circ}$ step trials for identification/training (mean and variance of steering angle calculated from these 14 trials; note that most likely the estimated variance does not represent the "true" variance of the human completely, only an approximation, see problems in representing the variance on validation data); 14 $-120^{\circ}$ step trials for validation (mean and variance calculated from these 14 trials)
		\item data preparation and cutting: cut data for the 14 participants, each trial individually according to threshold for steering angle velocity, start time = last time when absolute velocity under threshold, end time = first time when velocity under threshold again, average end time = time horizon N (determined by those trials which fulfill the end condition), trial is dropped if even relaxed start time condition cannot be fulfilled (movements starts not from an equilibrium condition and hence, it is assumed that not a single repetition is performed), calculate mean and variance curves for the identification and validation steps, reaction as well as premovement times are compensated via this procedure, indices of start time could be used to evaluate premovement times (due to look ahead in experimental design reaction times should be nearly negligible); in general, reaction and premovement times should be dropped for ISOC evaluation always; calculate starting mean values for velocity, torque and filter from angle and angle speed
		\item data analysis: discuss kinematic features (single-peaked, bell-shaped, nearly symmetric velocity profiles + max/mean ratio in range of reaching tasks + most profiles left-skewed, some nearly symmetric) \cite{Engelbrecht.2001} and thus the conclusion to model task as point-to-point movement and the hypothesis that control-dependent noise processes are crucial to describe the mean and variance curves
	\end{itemize}
}


\section{Results} \label{sec:results}
\begin{figure}[t]
	\centering
	\includegraphics[width=3in]{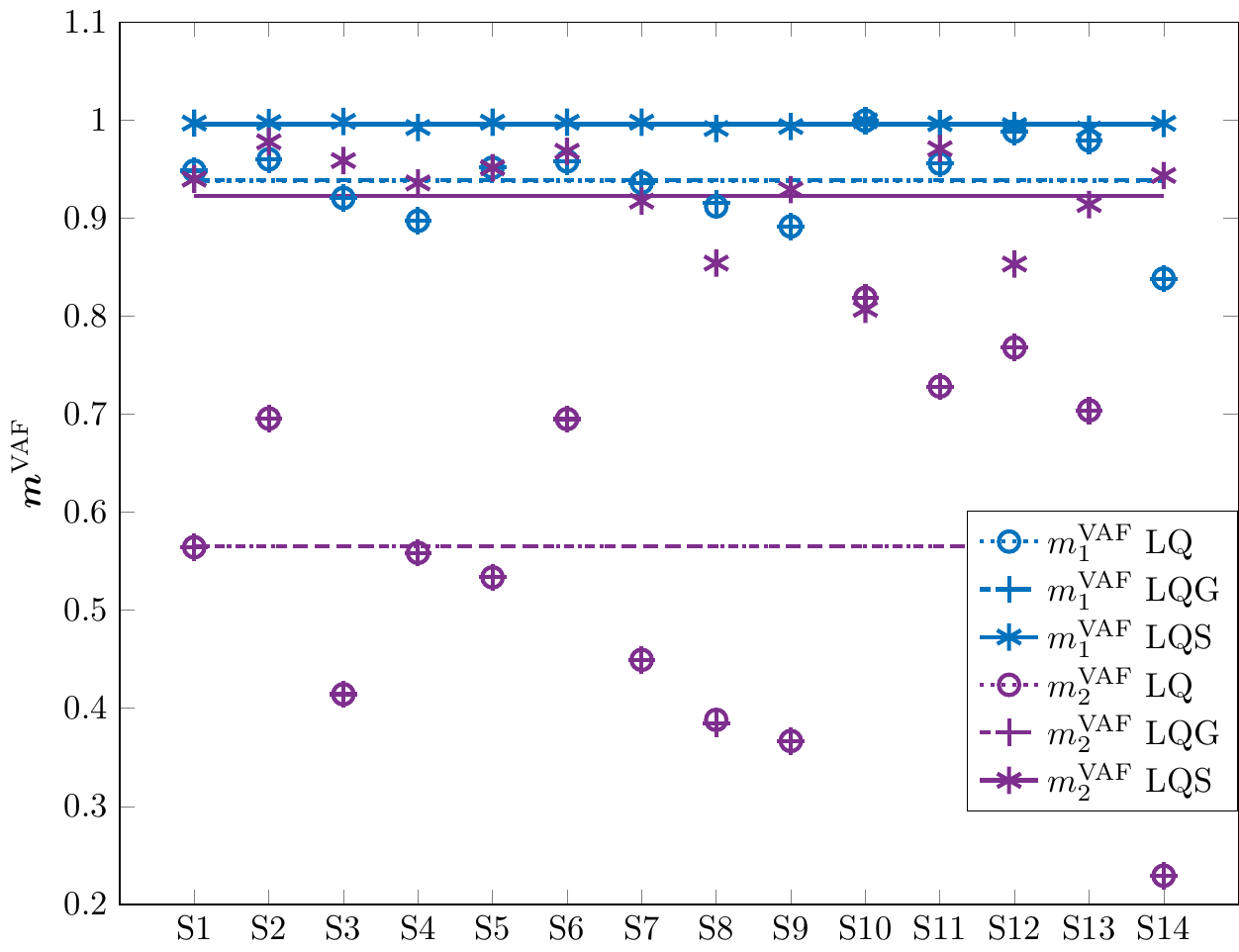}
	\caption{VAF values for the steering angle mean (blue) and steering angle velocity mean (purple) achieved with the identified LQS, LQG and LQ model on the training data set. Lines correspond to the average value over all subjects.}
	\label{fig:VAFmeanTrain}
\end{figure}
\begin{figure}[t]
	\centering
	\includegraphics[width=3in]{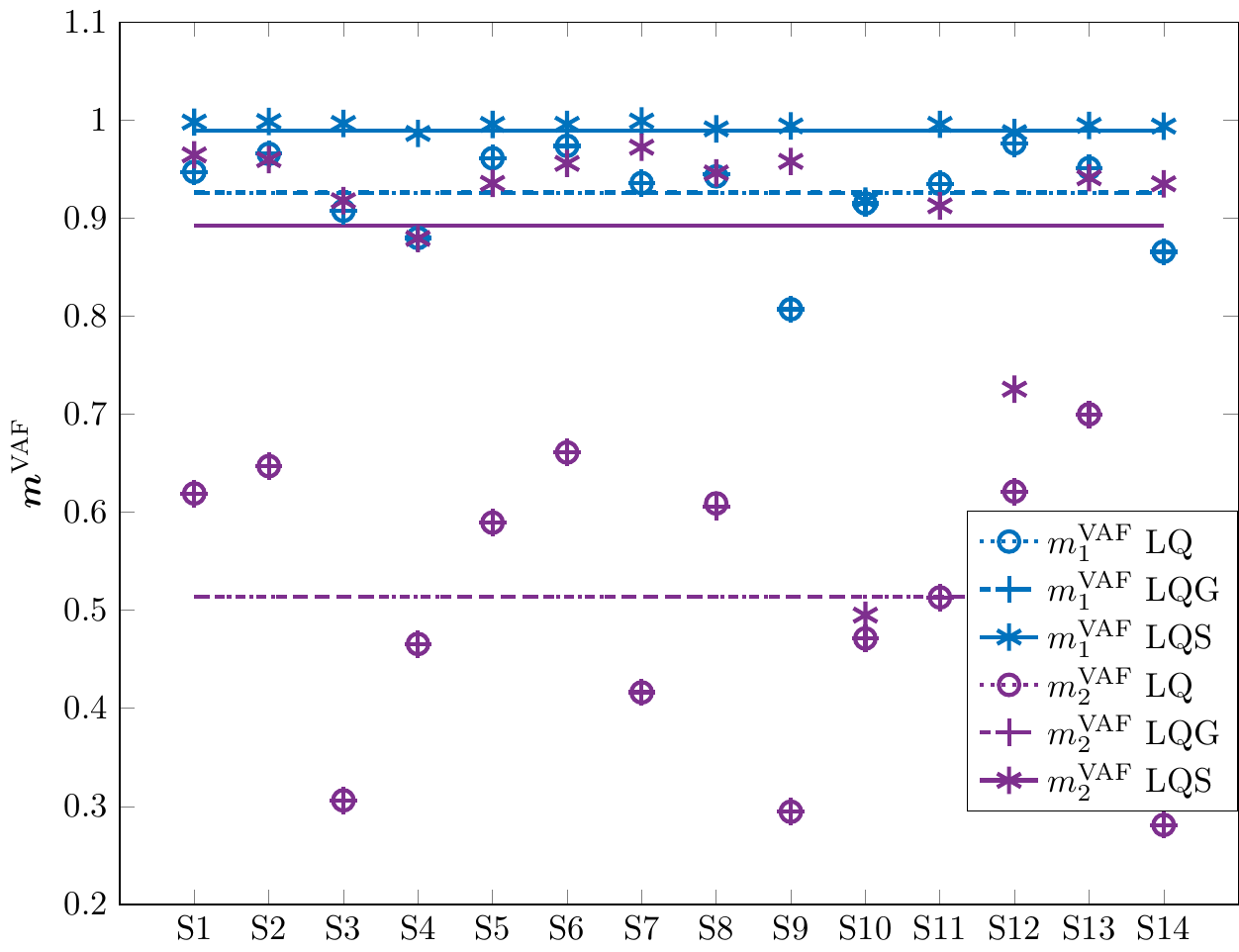}
	\caption{VAF values for the steering angle mean (blue) and steering angle velocity mean (purple) achieved with the identified LQS, LQG and LQ model on the validation data set. Lines correspond to the average value over all subjects.}
	\label{fig:VAFmeanTest}
\end{figure}
\begin{figure}[t]
	\centering
	\includegraphics[width=2.5in]{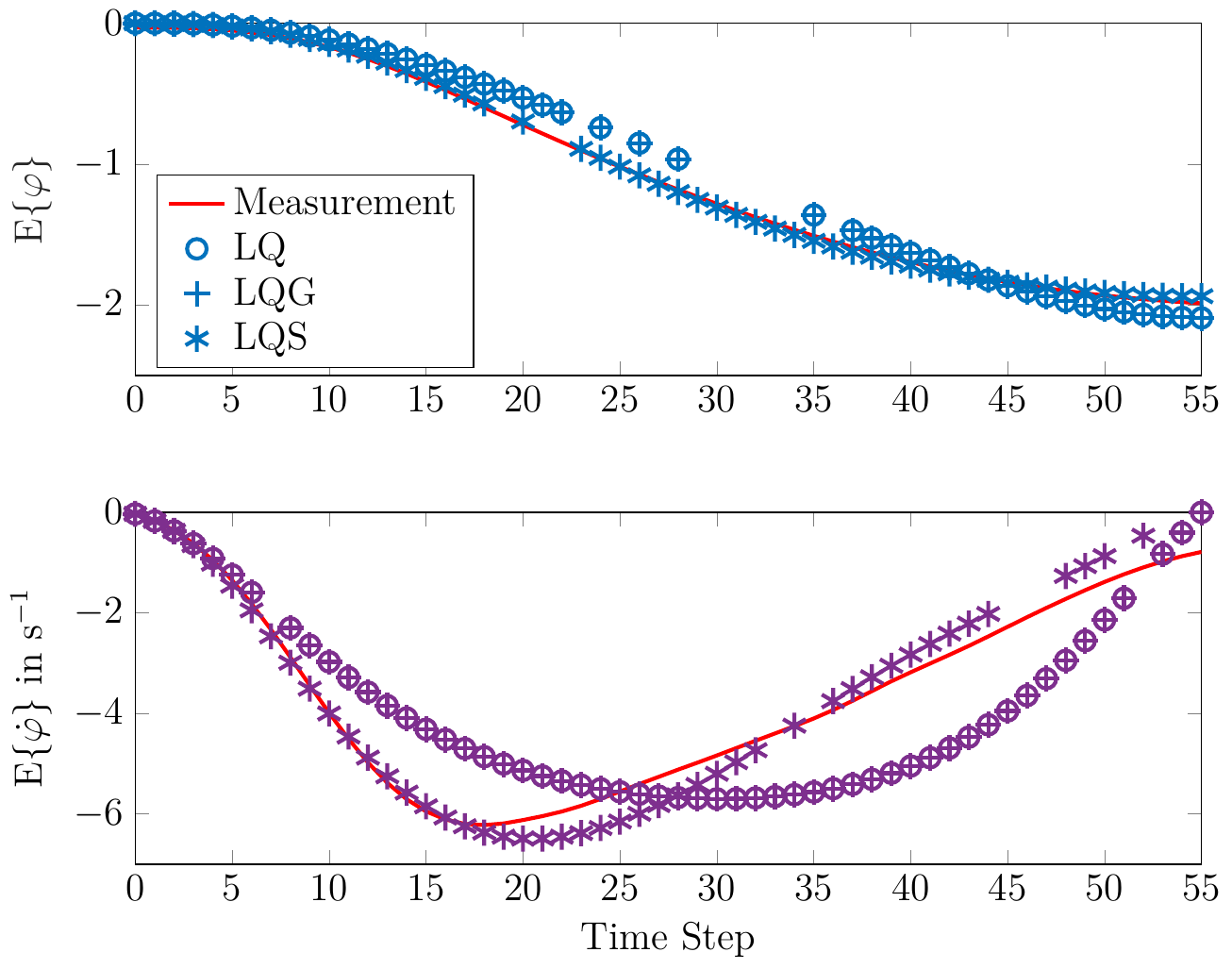}
	\caption{Predictions of all three models compared to human measurement data on validation data for subject 2.}
	\label{fig:meanComparisonS2}
\end{figure}

As motivated before, the focus of our model validation and comparison particularly lies on human average behavior. Hence, we first present detailed results of the LQS, LQG and LQ model in describing the human average steering behavior on training and validation data. Fig.~\ref{fig:VAFmeanTrain} and Fig.~\ref{fig:VAFmeanTest} show the VAF values of all 14 subjects as well as the average over all subjects for $\E{\varphi}$ and $\E{\dot{\varphi}}$ on training and validation data, respectively. The LQS model outperforms the LQG and LQ model for every subject and in average on training and validation data, regarding $\E{\varphi}$ and $\E{\dot{\varphi}}$. On validation data, the LQS model predictions show VAF values of $0.989$ ($\E{\varphi}$) and $0.893$ ($\E{\dot{\varphi}}$) in average compared to $0.926$ and $0.514$ in case of the LQG and LQ model. When considering average behavior, the LQG and LQ model show identical predictions and VAF values as expected since the additive noise processes do not influence the predicted steering angle (velocity) mean (see Corollary~\ref{corollary:meanLQS}). We perform an analysis of variance (ANOVA) to underline that the improvement of the LQS model is statistical significant by achieving $p$-values of $1.6 \cdot 10^{-4}$ and $2.5 \cdot 10^{-9}$ for $\E{\varphi}$ and $\E{\dot{\varphi}}$, respectively. On training data, very similar VAF and $p$-values result. 

In Fig.~\ref{fig:meanComparisonS2} the predictions of $\E{\varphi}$ and $\E{\dot{\varphi}}$ of the three models are depicted together with the human measurement for subject 2 to show the better performance of the LQS model qualitatively as well. Regarding the steering angle velocity, the LQS model especially describes the bell-shape and the maximum of the measurement curve better. 

\begin{figure}[t]
	\centering
	\includegraphics[width=2.5in]{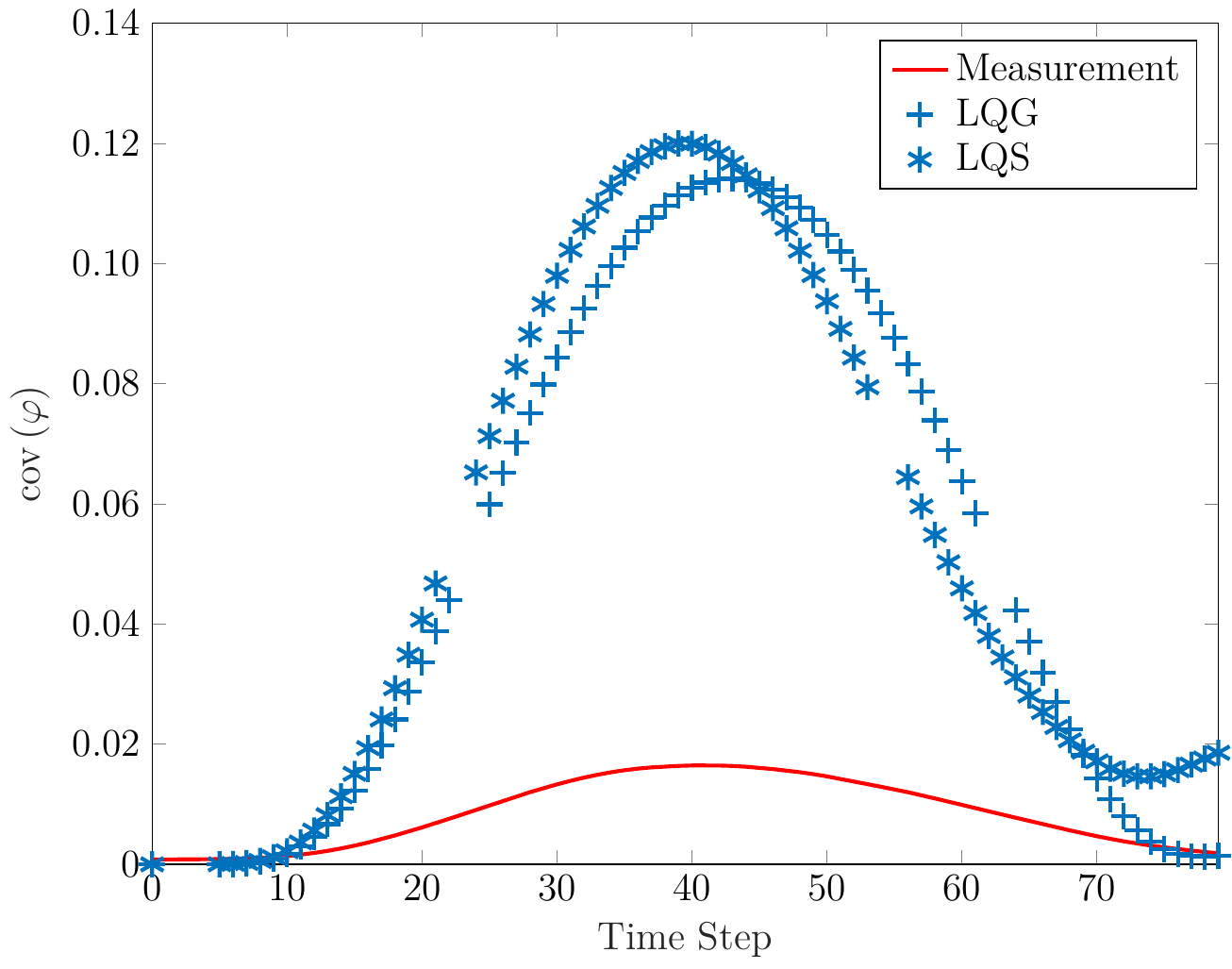}
	\caption{Predictions of the LQG and LQS model compared to human measurement data on validation data for subject 1.}
	\label{fig:varianceComparisonS1}
\end{figure}

\tdd{give complete VAF value ranges regarding variance in the following, explain/comment on actual variance of the movements in discussion section, add figure for example-wise results of subject 2, only steering angle, underline reason for that in the discussion section}
Concerning the steering angle variance, the human measurement data can only be described by the LQG and LQS model via overfitting. Both models achieve on the training data set VAF values from $0.451$ to $0.996$ with $0.880$ in average (LQG) and from $0.616$ to $0.985$ with $0.904$ in average (LQS), with no statistical difference. On validation data, the VAF values range from $-104.2$ to $0.895$ with $-7.288$ in average (LQG) and from $-107.1$ to $0.894$ with $-7.648$ in average (LQS). In Fig.~\ref{fig:varianceComparisonS1}, the predictions of the LQG and LQS model for subject 1 underline this poor performance in describing the measured steering angle variance qualitatively. 

\tdd{
	\begin{itemize}
		\item up to and including half of page 5
		\item VAF comparison regarding mean of steering angle and steering angle velocity for all three models (figure for training and validation data), results of ANOVA analysis for both quantities to show statistical significance (for training and validation data)
		\item plots for one example subject, i.e. S2, (2x1 plot), only validation data, ground truth/model predictions 
		\item VAF comparison regarding variance of steering angle and LQG, LQS model (name results for training data, name/explain differences on validation data), results of ANOVA analysis (for training data)
		\item show achieved scaling parameters for control-dependent noise (plot VAF mean steering angle and steering angle velocity over scaling parameter control-dependent noise, integrate box plot for parameters without outliers at averaged VAF, outliers: problematic subject 10 (see also achieved VAF values, training and validation data) due to atypical behavior (see kinematic feature analysis), two local optima achieved)
	\end{itemize}
}


\section{Discussion} \label{sec:discussion}
One of our main research objectives in this work is to analyze and prove the positive influence of the signal-dependent noise processes of the LQS model on describing human average behavior. Therefore, we begin the discussion of the results presented in Section~\ref{sec:results} by formulating the following hypothesis.
\begin{hypothesis} \label{hypothesis:averageSteeringBehavior}
	The LQS model describes human average steering behavior represented by the mean of steering angle $\E{\varphi}$ and steering angle velocity $\E{\dot{\varphi}}$ significantly better than a LQG or deterministic LQ model.	
\end{hypothesis}
With Fig.~\ref{fig:VAFmeanTest} and the ANOVA results, Hypothesis~\ref{hypothesis:averageSteeringBehavior} can be accepted since the LQS model outperforms the others for significance levels well below $1\,\%$. Remarkably, the LQS model is able to describe the velocity profiles (see Fig.~\ref{fig:velocityProfiles}) better since the average VAF value for $\E{\dot{\varphi}}$ is nearly $75\,\%$ higher in case of the LQS model. Furthermore, the qualitative shape of the profiles matches better (see Fig.~\ref{fig:meanComparisonS2}). However, $\E{\dot{\varphi}}$ is predicted slightly worse than $\E{\varphi}$ which can be concluded from the VAF values in Fig.~\ref{fig:VAFmeanTrain} and Fig.~\ref{fig:VAFmeanTest} as well as from example subject 2 (see Fig.~\ref{fig:meanComparisonS2}). The main reason for this is that the sections with constant reference $\varphi_{\text{r}}$ were too short. Consequently, we set the threshold for determining the movement duration ($\dot{\varphi}<\SI{0.1}{s^{-1}}$) relatively high to separate a single movement from a movement to a new target steering angle. Therefore, several velocity profiles do not show a convergence to zero at the movement ending which however the LQS model tends to predict (see Fig.~\ref{fig:meanComparisonS2}). Furthermore, subject 10 seems to have problems in following these tough reference changes resulting in its outlier behavior (see Fig.~\ref{fig:velocityProfiles}) which none of the models can describe (VAF values $<0.5$, see Fig.~\ref{fig:VAFmeanTest}). Hence, with a reworked study design regarding the target steering angle we expect similar results for the predictions of $\E{\varphi}$ and $\E{\dot{\varphi}}$ in case of the LQS model as well as no outlier behavior of participants.

\begin{figure}[t]
	\centering
	\includegraphics[width=2.5in]{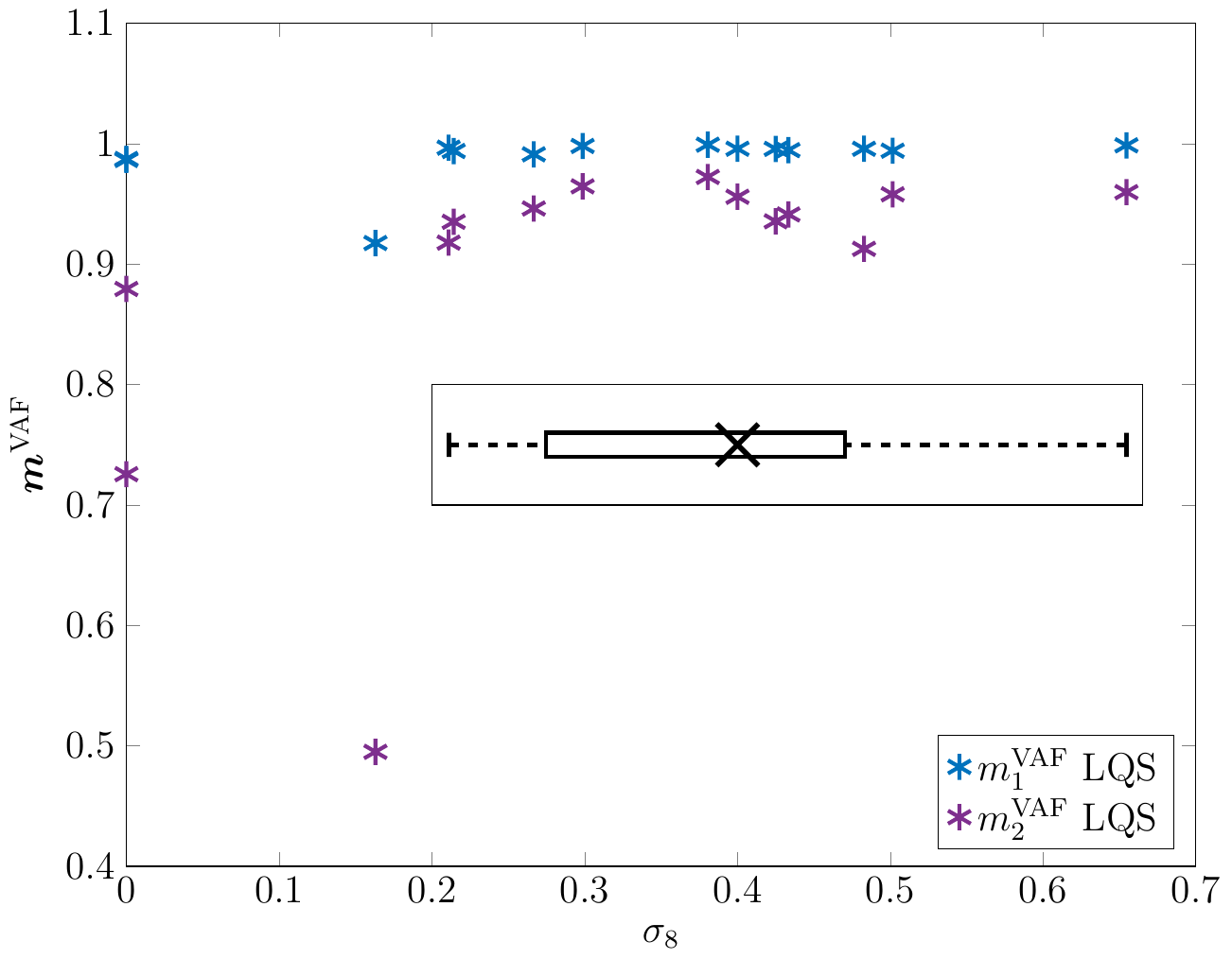}
	\caption{VAF values achieved with the LQS model on validation data over control-dependent noise scaling parameter $\sigma_8$. Box plot for $\sigma_8$ is shown additionally where outliers at $\sigma_8 \approx 0$ and $\sigma_8 \approx 0.16$ are dropped.}
	\label{fig:VAFoverControlNoiseScaling}
\end{figure}
The significant better modeling performance of the LQS model compared to the reduced versions of it without signal-dependent noise processes is already a strong indication that especially these parts of the model have a positive influence on describing human average behavior. In order to provide a further evidence for this claim, Fig.~\ref{fig:VAFoverControlNoiseScaling} depicts an analysis of the scaling parameter of the control-dependent noise process which we assume to have the biggest positive influence according to findings in literature (see \cite{Harris.1998}). In Fig.~\ref{fig:VAFoverControlNoiseScaling}, the VAF values for $\E{\varphi}$ and $\E{\dot{\varphi}}$ achieved on validation data with the identified LQS model are plotted over $\sigma_8$. Now, when the data points at $\sigma_8 \approx 0.16$ (subject 10) and at $\sigma_8 \approx 0$ (ISOC algorithm converged into a local optimum) are considered as outliers, for $\sigma_8$ the box plot in Fig.~\ref{fig:VAFoverControlNoiseScaling} results. The box plot highlights that the high VAF values ($0.996$ and $0.945$ in average for $\E{\varphi}$ and $\E{\dot{\varphi}}$) of the remaining subjects are achieved with $\sigma_8$-values greater than $0.2$. Hence, our data-driven analysis confirms findings in literature and theoretically expected by Corollary~\ref{corollary:meanLQS} that the signal-dependent and particularly the control-dependent noise process of the LQS model improves the performance of modeling human average behavior.

Concerning the description of human variability patterns with different SOC models, we cannot draw conclusions from our results. The main reason is that 14 repetitions of the same movement are (most likely) not enough to approximate the real human variability pattern. Hence, it cannot be assumed that the assumption on the estimate $\hat{\bm{\Omega}}_t^{\bm{x}^{\ast}}$ in Problem~\ref{problem:ISOC} is fulfilled and the parameters determined on training data show poor generalization related to covariance predictions. Fig.~\ref{fig:varianceComparisonS1} yields a possible reason for this behavior. Although the measured variance of the steering angle in case of the shown $-\frac{2}{3}\pi$-step in the reference (validation data) shows reasonable characteristics, i.e. single-peaked curve with a maximum at the maximum velocity, the maximum value is with $0.0165$ around one order of magnitude smaller than that on the training data ($0.139$). Since different but in the order of magnitude identical variance curves for the $+\frac{2}{3}\pi$-step and $-\frac{2}{3}\pi$-step would be expected, the estimated variance seems to be not representative for the real human variability. In future work, we will increase the number of single movements observed from subjects to overcome this problem. Since this poor variance description can already be seen at the steering angle, we dropped the velocity, which is only numerically computed from it, in the covariance fitting and analysis.

Although we look at a very specific application field (human steering behavior) at a first glance, the analysis of the kinematic features highlights that steering in our simplified driving task is indeed comparable to point-to-point human hand reaching (see Section~\ref{subsec:driving_experiment}). Hence, we strengthen on one side the proposition done by \cite{Kolekar.2018} that steering is comparable to reaching and can make on the other a step towards a data-driven validation and analysis of SOC models for general goal-directed human movements. Thus, we are strongly convinced to be able to draw similar conclusions when we apply our ISOC algorithm to data of human hand movements: significant better performance of the LQS model in describing human average behavior quantitatively as well as qualitatively and scaling parameters of the control-dependent noise process as the main reason for this modeling performance.


\tdd{
\begin{itemize}
	\item up to page 6
	\item hypothesis 1: LQS best model to describe human average steering behavior (steering angle and steering angle velocity)
	\item conclusion 1: hypothesis 1 accepted and via box plot of scaling parameter control-dependent noise this can be traced back to the control-dependent noise process (see \cite{Harris.1998}); points for future work: longer single movement phases (see example plots for one subject, GT not in rest at the end + problematic subject 10) + better sensor setup (better numerical calculation of velocity) + improved algorithm (local optima)
	\item conclusion 2: no conclusion regarding variance possible; future work: more data (variance analysis, this also the reason for setting the weighting vectors appropriately, since this problem is already present for steering angle, we dropped the steering angle velocity as directly derived signal in the variance analysis and variance fitting, with more data and a better sensor setup we want to overcome this as well)
	\item conclusion 3: steering as reaching (see kinematic features and the model validation in this task although models were developed/investigated mainly on point-to-point or via-point movements), majority of velocity profiles left-skewed reasonable since faster movements preferred due to definition of task (see \cite{Engelbrecht.2001}), analysis strengths and deepens discussion in \cite{Kolekar.2018}
	\item conclusion 4: LQ model not so worse overall, together with conclusion 3 this explains relative good performance when it is extended to a differential game model for cooperative driving between two humans \cite{Inga.2019c}
\end{itemize}
}


\section{Conclusion} \label{sec:conclusion}
In this paper, we use our previously developed Inverse Stochastic Optimal Control (ISOC) algorithm solving the identification problem of the linear-quadratic (LQ) sensorimotor (LQS) model to compare its capability in describing human steering behavior in a simplified driving task to a LQ Gaussian (LQG) and deterministic LQ model. In this context, an important focus lies on human average behavior. The final results show that the LQS model is able to predict the mean of steering angle and steering angle velocity with Variance Accounted For (VAF) values of $0.989$ and $0.893$ on a validation data set and averaged over 14 subjects. The improvement compared to the results of the LQG and LQ model is statistical significant with significance levels well below $1\,\%$. Based on a parameter analysis, this significant better performance can be traced back to the control-dependent noise process characteristic for the LQS model. Since the movements in our simplified driving task show similar kinematic features as point-to-point human hand movements, this is a strong indication that the signal-dependent and in particular the control-dependent noise process of the LQS model is one of the crucial model parts to describe human average behavior. For human-machine systems, this better performance of the LQS model in characterizing human average behavior can be used for more accurate predictions of human movements in shared workspace applications or to facilitate a more intuitive support of humans in a trajectory tracking task.
\tdd{comment on possible improvements of using identified LQS models for human-machine cooperation}


\td{Update version on arxiv}

\bibliography{References2}             

\begin{thebibliography}{18}
\providecommand{\natexlab}[1]{#1}
\providecommand{\url}[1]{\texttt{#1}}
\providecommand{\urlprefix}{URL }
\expandafter\ifx\csname urlstyle\endcsname\relax
  \providecommand{\doi}[1]{doi:\discretionary{}{}{}#1}\else
  \providecommand{\doi}{doi:\discretionary{}{}{}\begingroup
  \urlstyle{rm}\Url}\fi

\bibitem[{Engelbrecht(2001)}]{Engelbrecht.2001}
Engelbrecht, S. (2001).
\newblock Minimum principles in motor control.
\newblock \emph{Journal of Mathematical Psychology}, 45, 497--542.

\bibitem[{Flash and Hogan(1985)}]{Flash.1985}
Flash, T. and Hogan, N. (1985).
\newblock The coordination of arm movements: An experimentally confirmed
  mathematical model.
\newblock \emph{The Journal of Neuroscience}, 5(7), 1688--1703.

\bibitem[{Franklin and Wolpert(2011)}]{Franklin.2011}
Franklin, D.W. and Wolpert, D.M. (2011).
\newblock Computational mechanisms of sensorimotor control.
\newblock \emph{Neuron}, 72(3), 425--442.

\bibitem[{Gallivan et~al.(2018)Gallivan, Chapman, Wolpert, and
  Flanagan}]{Gallivan.2018}
Gallivan, J.P., Chapman, C.S., Wolpert, D.M., and Flanagan, J.R. (2018).
\newblock Decision-making in sensorimotor control.
\newblock \emph{Nature Reviews. Neuroscience}, 19(9), 519--534.

\bibitem[{Harris and Wolpert(1998)}]{Harris.1998}
Harris, C.M. and Wolpert, D.M. (1998).
\newblock Signal-dependent noise determines motor planning.
\newblock \emph{Nature}, 394, 780--784.

\bibitem[{Karg et~al.(2021)Karg, Kalb, Bengler, and Hohmann}]{Karg.2021}
Karg, P., Kalb, L., Bengler, K., and Hohmann, S. (2021).
\newblock Passing the baton between autopilot and driver.
\newblock \emph{ATZ Worldwide}, 123, 68--72.

\bibitem[{Karg et~al.(2022)Karg, Stoll, Rothfu{\ss}, and Hohmann}]{Karg.2022}
Karg, P., Stoll, S., Rothfu{\ss}, S., and Hohmann, S. (2022).
\newblock Inverse stochastic optimal control for linear-quadratic gaussian and
  linear-quadratic sensorimotor control models.
\newblock \emph{IEEE Conference on Decision and Control}, 2801--2808.

\bibitem[{Kolekar et~al.(2018)Kolekar, Mugge, and Abbink}]{Kolekar.2018}
Kolekar, S., Mugge, W., and Abbink, D. (2018).
\newblock Modeling intradriver steering variability based on sensorimotor
  control theories.
\newblock \emph{IEEE Transactions on Human-Machine Systems}, 48(3), 291--303.

\bibitem[{MacAdam(1981)}]{MacAdam.1981}
MacAdam, C.C. (1981).
\newblock Application of an optimal preview control for simulation of
  closed-loop automobile driving.
\newblock \emph{IEEE Transactions on Systems, Man, and Cybernetics}, 11(6),
  393--399.

\bibitem[{Markkula(2014)}]{Markkula.2014}
Markkula, G. (2014).
\newblock Modeling driver control behavior in both routine and near-accident
  driving.
\newblock \emph{Human Factors and Ergonomics Society Annual Meeting}, 58(1),
  879--883.

\bibitem[{Nash and Cole(2015)}]{Nash.2015}
Nash, C.J. and Cole, D.J. (2015).
\newblock Development of a novel model of driver-vehicle steering control
  incorporating sensory dynamics.
\newblock \emph{Dynamics of Vehicles on Roads and Tracks}, 72--81.

\bibitem[{Ross and Bagnell(2010)}]{Ross.2010}
Ross, S. and Bagnell, J.A. (2010).
\newblock Efficient reductions for imitation learning.
\newblock \emph{International Conference on Artificial Intelligence and
  Statistics}, 661--668.

\bibitem[{Sentouh et~al.(2009)Sentouh, Chevrel, Mars, and
  Claveau}]{Sentouh.2009}
Sentouh, C., Chevrel, P., Mars, F., and Claveau, F. (2009).
\newblock A sensorimotor driver model for steering control.
\newblock \emph{IEEE International Conference on Systems, Man, and
  Cybernetics}, 2462--2467.

\bibitem[{Todorov(2004)}]{Todorov.2004}
Todorov, E. (2004).
\newblock Optimality principles in sensorimotor control.
\newblock \emph{Nature Neuroscience}, 7(9), 907--915.

\bibitem[{Todorov(2005)}]{Todorov.2005}
Todorov, E. (2005).
\newblock Stochastic optimal control and estimation methods adapted to the
  noise characteristics of the sensorimotor system.
\newblock \emph{Neural Computation}, 17, 1084--1108.

\bibitem[{Todorov and Jordan(2002)}]{Todorov.2002}
Todorov, E. and Jordan, M.I. (2002).
\newblock Optimal feedback control as a theory of motor coordination.
\newblock \emph{Nature Neuroscience}, 5(11), 1226--1235.

\bibitem[{Uno et~al.(1989)Uno, Kawato, and Suzuki}]{Uno.1989}
Uno, Y., Kawato, M., and Suzuki, R. (1989).
\newblock Formation and control of optimal trajectory in human multijoint arm
  movement.
\newblock \emph{Biological Cybernetics}, 61, 89--101.

\bibitem[{Winter(1990)}]{Winter.1990}
Winter, D.A. (1990).
\newblock \emph{Biomechanics and Motor Control of Human Movement}.
\newblock Wiley, 2nd edition.

\end{thebibliography}
                                                   







\appendix

\end{document}